\documentclass[12pt, english, a4paper]{article}
\usepackage[latin9]{inputenc}
\usepackage[T1]{fontenc}
\usepackage{graphicx}
\usepackage{amsthm}
\usepackage{amsfonts}
\usepackage{circledsteps}
\usepackage{mathtools}
\allowdisplaybreaks
\usepackage{amssymb}
\usepackage{makecell}
\usepackage{comment}
\usepackage{booktabs}
\usepackage{authblk}
\usepackage{cases}
\usepackage{hyperref}
\usepackage{mathrsfs}
\usepackage{fullpage}
\usepackage{times}
\usepackage{soul}
\usepackage{booktabs}
\usepackage{subcaption}
\usepackage[dvipsnames]{xcolor}
\usepackage{hyperref}
\usepackage{tikz,tikz-cd}

\usepackage[colorinlistoftodos, prependcaption, textsize=tiny]{todonotes}
\usepackage{xcolor} 

\usetikzlibrary{decorations.text}

\usetikzlibrary{arrows.meta,positioning}
\usepackage{extarrows}
\usepackage{float}
\usepackage{ytableau}
\usepackage{xcolor}
\usepackage{circledsteps}
\usepackage{multirow}
\usepackage{bm}
\definecolor{pcdcolor}{RGB}{255, 230, 230}

\theoremstyle{plain}
\newtheorem{theorem}{Theorem}

\newtheorem{lemma}[theorem]{Lemma}

\theoremstyle{definition}
\newtheorem{definition}[theorem]{Definition}
\newtheorem{example}{Example}

\theoremstyle{remark}

\usepackage{authblk}

\title{\textbf{{A Combinatorial Model for the Faces of Graphical Threshold Arrangements}}}

\author[]{{\small Shubhanshu Prasad}\thanks{\tt{shubhanshup@iisc.ac.in}}}
\affil[]{{\small Department of Mathematics, Indian Institute of Science, Bengaluru 560012, India}}

\date{}

\begin{document}

\maketitle

\begin{abstract}
For a graph $G$ on the vertex set $[n]$, the $G$-threshold arrangement $\mathcal{T}(G)$ is the arrangement in $\mathbb{R}^n$ consisting of the hyperplanes $x_i+x_j=0$ for $ij\in E(G)$. For $G=\mathbb{K}_n$, this is the usual threshold arrangement, whose regions are known to be in bijection with labelled threshold graphs. We generalize this correspondence to arbitrary graphs $G$ by introducing $G$-threshold graphs, constructed from signed sequences, and showing that they are in bijection with the regions of $\mathcal{T}(G)$. We also describe all faces of $\mathcal{T}(G)$ using a combinatorial model of colored signed sequences and colored $G$-threshold graphs, and show that the dimension of a face is the number of bipartite connected components of the graph formed by the edges within the color classes. Using this model, we obtain formulas for the number of faces of each dimension of $\mathcal{T}(\mathbb{K}_n)$ in terms of compositions, Fubini numbers, and the numbers of labelled threshold graphs, and for $\mathcal{T}(\mathbb{K}_{m,n})$, where the enumeration reduces to counting acyclic orientations of the join of a complete graph and a complete bipartite graph.
\end{abstract}

\textbf{\small{Keyword:}}{\small{ hyperplane arrangements, threshold graphs, faces and regions, signed sequences, acyclic orientations, chromatic polynomial} }{\let\thefootnote\relax}

\section{Introduction}

An affine hyperplane $H$ in $\mathbb{R}^n$ is the set of all points $\mathbf{x}\in\mathbb{R}^n$ satisfying $\langle \mathbf{a},\mathbf{x}\rangle=c$, where $\mathbf{a}\in\mathbb{R}^n\setminus\{0\}$ and $c\in\mathbb{R}$ are fixed. A \emph{hyperplane arrangement} $\mathcal{A}$ is a finite collection of affine hyperplanes in $\mathbb{R}^n$ (see \cite{Stanley2007} for background). Let $\mathcal{A}=\{H_i\}_{i\in[m]}$, where the defining equation of $H_i$ is $\langle \mathbf{a}_i,\mathbf{x}\rangle=c_i$, with $\mathbf{a}_i\in\mathbb{R}^n$ and $c_i\in\mathbb{R}$. A \emph{face} of $\mathcal{A}$ is a nonempty set of the form
\[
F_{\varepsilon}=
\left\{
\mathbf{x}\in\mathbb{R}^n:
\operatorname{sgn}(\langle \mathbf{a}_i,\mathbf{x}\rangle-c_i)
=\varepsilon_i \text{ for all } i\in[m]
\right\},
\]
where $\varepsilon\in\{-,0,+\}^m$. A face corresponding to a sign vector $\varepsilon\in\{-,+\}^m$ does not intersect any hyperplane of the arrangement and is called a \emph{region}.

The enumeration of regions and faces of hyperplane arrangements has been studied extensively. Many classical combinatorial objects arise naturally in this setting. For instance, the regions of the braid arrangement are in bijection with permutations, while the regions of the Shi arrangement are in bijection with parking functions \cite{Shi1986,Athanasiadis1999}. The finite field method provides a powerful tool for enumerating regions and faces of hyperplane arrangements \cite{Athanasiadis1996,Stanley2007}. In addition to enumerative results, bijective approaches are useful for understanding the combinatorial structure underlying these geometric objects. In particular, they provide explicit correspondences between the regions or faces of an arrangement and suitable combinatorial models.

In this paper, we study the \emph{threshold arrangement} $\mathcal{T}(\mathbb{K}_n)$, defined by
\[
\mathcal{T}(\mathbb{K}_n)=
\left\{
\left\{\mathbf{x}\in\mathbb{R}^n:x_i+x_j=0\right\}
\mid 1\leq i<j\leq n
\right\}.
\]
A graph is called a \emph{threshold graph} if its vertices can be added one at a time so that each new vertex is either adjacent to all the earlier vertices or to none of them \cite{MahadevPeled1995}. A threshold graph with vertex set $[n]$ is called a \emph{labelled threshold graph}. The recursive construction of threshold graphs gives a natural representation by sign sequences: a positive sign corresponds to adding a vertex adjacent to all earlier vertices, while a negative sign corresponds to adding an isolated vertex. It is known that the regions of $\mathcal{T}(\mathbb{K}_n)$ are in bijection with labelled threshold graphs on $[n]$. This connection has been used to study the characteristic polynomial of the threshold arrangement and the enumeration of labelled threshold graphs; see, for example, \cite{DeshpandeMenonSingh2021}. The correspondence between the regions of $\mathcal{T}(\mathbb{K}_n)$ and labelled threshold graphs describes only the $n$-dimensional faces of the arrangement. In this paper, we provide a combinatorial description of all faces of $\mathcal{T}(\mathbb{K}_n)$ and enumerate them according to dimension, in the more general setting described next.

For a graph $G$ with vertex set $[n]$, we define the \emph{$G$-threshold arrangement} by
\[
\mathcal{T}(G)=
\left\{
\left\{\mathbf{x}\in\mathbb{R}^n:x_i+x_j=0\right\}
\mid ij\in E(G)
\right\}.
\]
Thus, $\mathcal{T}(\mathbb{K}_n)$ is the usual threshold arrangement, and the correspondence between its regions and labelled threshold graphs raises a natural question: what happens for an arbitrary graph $G$? When $G$ is bipartite, the answer is classical. Changing the signs of the coordinates corresponding to one part of a bipartition transforms $\mathcal{T}(G)$ into the graphical arrangement of $G$, whose regions are in bijection with the acyclic orientations of $G$ \cite{Stanley2007}. For graphs with odd cycles this is no longer possible, and a new feature appears: on a face, each connected component of the graph formed by the edges $ij$ with $x_i+x_j=0$ leaves one free parameter if it is bipartite, and forces all its coordinates to be zero if it contains an odd cycle. The threshold arrangement $\mathcal{T}(\mathbb{K}_n)$ is the extreme case in which every pair of vertices is an edge. This suggests studying $\mathcal{T}(G)$ for an arbitrary graph $G$, and describing its regions and faces in terms of combinatorial objects that generalize labelled threshold graphs and take the zero coordinates into account.

Our first main contribution is a combinatorial model for the regions, called \emph{$G$-threshold graphs}. A $G$-threshold graph $T_G(\pi,\epsilon)$ is obtained from a signed sequence $(\pi,\epsilon)$ by adding the vertices one at a time and joining a vertex of sign $+$ to those earlier vertices to which it is adjacent in $G$. We show that two signed sequences give the same $G$-threshold graph if and only if they are $G$-equivalent, so that $\mathcal{S}_n(G)$, the set of equivalence classes, is in bijection with $\mathrm{Thr}(G)$. We then show that $\mathcal{S}_n(G)$ is also in bijection with the set $\mathcal{R}(\mathcal{T}(G))$ of regions of $\mathcal{T}(G)$ (Theorem~\ref{thm:regions_threshold_graphs}). This generalizes the correspondence between labelled threshold graphs and the regions of $\mathcal{T}(\mathbb{K}_n)$.

To describe all faces, and not only the regions, we must allow some of the equalities $x_i+x_j=0$ to hold. We therefore introduce \emph{colored signed sequences} $(\pi,\epsilon,\lambda)$ and the associated \emph{colored $G$-threshold graphs}. In a colored signed sequence, the composition $\lambda$ groups the vertices into color classes, which are the sets of vertices whose coordinates have the same absolute value on the face, and a vertex may have the sign $0$. The colored $G$-threshold graph records the edges inside the color classes, the pairs inside a color class whose endpoints have the same nonzero sign, and the edges between the color classes. We show that the set $\mathcal{C}\mathrm{Thr}(G)$ of colored $G$-threshold graphs is in bijection with the set $\mathcal{F}(\mathcal{T}(G))$ of faces of $\mathcal{T}(G)$ (Theorem~\ref{thm:faces_threshold_graphs}). Furthermore, we show that the dimension of a face is determined by the number of bipartite connected components of the graph formed by the edges within the color classes of its corresponding colored $G$-threshold graph (Theorem~\ref{thm:dim_through_thresh_grph}).

Finally, we use this model to count the faces of $\mathcal{T}(G)$ of each dimension in two cases. For $G=\mathbb{K}_n$, we obtain a formula in terms of strong compositions $\lambda$, Fubini numbers $a(r)$, and the numbers $|\mathrm{Thr}(\mathbb{K}_r)|$ of labelled threshold graphs (Theorem~\ref{thm:faces_of_thr_arrang}). For $G=\mathbb{K}_{m,n}$, we show that, for a fixed set of color classes, the colored $G$-threshold graphs are in bijection with the acyclic orientations of the graph $H_{\mathcal{C}}$ on the color classes, which is the join of a complete graph and a complete bipartite graph. Counting these with chromatic polynomials and Stanley's theorem on acyclic orientations \cite{Stanley1973}, we obtain a formula for the number of faces of each dimension (Theorem~\ref{thm:faces_Kmn}).

The paper is organized as follows. In Section~2, we define $G$-threshold graphs and establish a bijection between them and the regions of $\mathcal{T}(G)$. The proof proceeds in two steps. First, we define an equivalence relation on signed sequences and show that two signed sequences give the same $G$-threshold graph if and only if they are equivalent. This gives a bijection between the equivalence classes of signed sequences and the set of $G$-threshold graphs. In the second step, we establish a bijection between these equivalence classes and the regions of $\mathcal{T}(G)$.

In Section~3, we extend the correspondence of Section~2 to all faces of $\mathcal{T}(G)$. We introduce colored signed sequences and colored $G$-threshold graphs to describe these faces. Following the approach of Section~2, we first define an equivalence relation on colored signed sequences and show that two colored signed sequences give the same colored $G$-threshold graph if and only if they are equivalent. We then establish a bijection between the equivalence classes of colored signed sequences and the faces of $\mathcal{T}(G)$. Finally, we determine the dimension of each face in terms of the bipartite connected components of the graph formed by the edges within the color classes of its corresponding colored $G$-threshold graph.

In the final section, we use this combinatorial model to derive formulas for the number of faces of each dimension of $\mathcal{T}(G)$ in two cases: when $G=\mathbb{K}_n$ and when $G=\mathbb{K}_{m,n}$.
\section{Bijection between regions of $G$-threshold arrangement and $G$-threshold graph}

We first recall the notion of a signed sequence. A \textit{signed sequence} is a pair $(\pi,\epsilon)$, where $\pi\in S_n$ is a permutation and $\epsilon=(\epsilon_1,\dots,\epsilon_n)\in\{+,-\}^n$, where the number $\pi(i)$ in the $i$th position of $\pi$ is assigned the sign $\epsilon_i$. We denote the set of all signed sequences by $\mathcal{S}_n$. For a signed sequence $(\pi,\epsilon)\in\mathcal S_n$, the \textit{sign} of a vertex $v\in[n]$ is $\epsilon_i$, where $i=\pi^{-1}(v)$ is the position of $v$ in $\pi$.

Recall that a threshold graph can be constructed from a signed sequence by adding the vertices in the given order, where a vertex with sign $-$ is added as an isolated vertex and a vertex with sign $+$ is joined to all previously added vertices. The following definition is a natural generalization of this construction, with the graph $G$ determining which edges can be added. In particular, when $G=\mathbb{K}_n$, we recover the usual threshold graphs.

\begin{definition}[$G$-threshold graph]\label{def:G_threshold_graph}
Let $G$ be a graph with vertex set $[n]$, and let $(\pi,\epsilon)$ be a signed sequence, where $\pi=\pi_1\pi_2\dots \pi_n$ is a permutation and $\epsilon=\epsilon_1\epsilon_2\dots\epsilon_n$ is a sequence of signs. We define $T_G(\pi,\epsilon)$ to be a \textit{labelled $G$-threshold graph} on vertex set $[n]$ obtained recursively as follows:
\begin{enumerate}
\item Start with the vertex $\pi_1$.
\item Suppose that the vertices $\pi_1,\dots,\pi_{i-1}$ have been added. Add the vertex $\pi_i$. If $\epsilon_i=-$, do not add any edges incident to $\pi_i$. If $\epsilon_i=+$, then, for each $k\in[i-1]$, add an edge between $\pi_k$ and $\pi_i$ whenever $\pi_k\pi_i\in E(G)$.
\end{enumerate}
The set of all the labelled $G$-threshold graphs is denoted by $\operatorname{Thr}(G)$.
\end{definition}

We illustrate Definition~\ref{def:G_threshold_graph} with an example. It shows how a $G$-threshold graph is built step by step, and that different signed sequences can give the same graph.

\begin{example}\label{ex:threshold_graph}
Let $n=4$ and let $G$ be the graph on $[4]$ with edge set $E(G)=\{12,\,13,\,23,\,34\}$, that is, a triangle on $\{1,2,3\}$ with a pendant edge $34$ (see Figure~\ref{fig:G_and_threshold}(a)).
Consider the signed sequences
\[
\begin{aligned}
(\pi,\epsilon)       &= (1234,\ (+,-,+,-)),\\
(\pi',\epsilon')     &= (1234,\ (-,-,+,-)),\\
(\pi'',\epsilon'')   &= (2134,\ (-,-,+,-)).
\end{aligned}
\]

\medskip

\noindent\textbf{Constructing the graph for $(\pi,\epsilon)$:} We start with the vertex $1$. Then we add $2$, which has sign $-$, so no edge is added. Next we add $3$, which has sign $+$, so we join it to its neighbors in $G$ among the earlier vertices $1,2$. Since $13,23\in E(G)$, we add $13$ and $23$. Finally, we add $4$, which has sign $-$, so no edge is added. The result is the graph $T_G(\pi,\epsilon)$ with edge set $\{13,\,23\}$ (see Figure~\ref{fig:G_and_threshold}(b)).

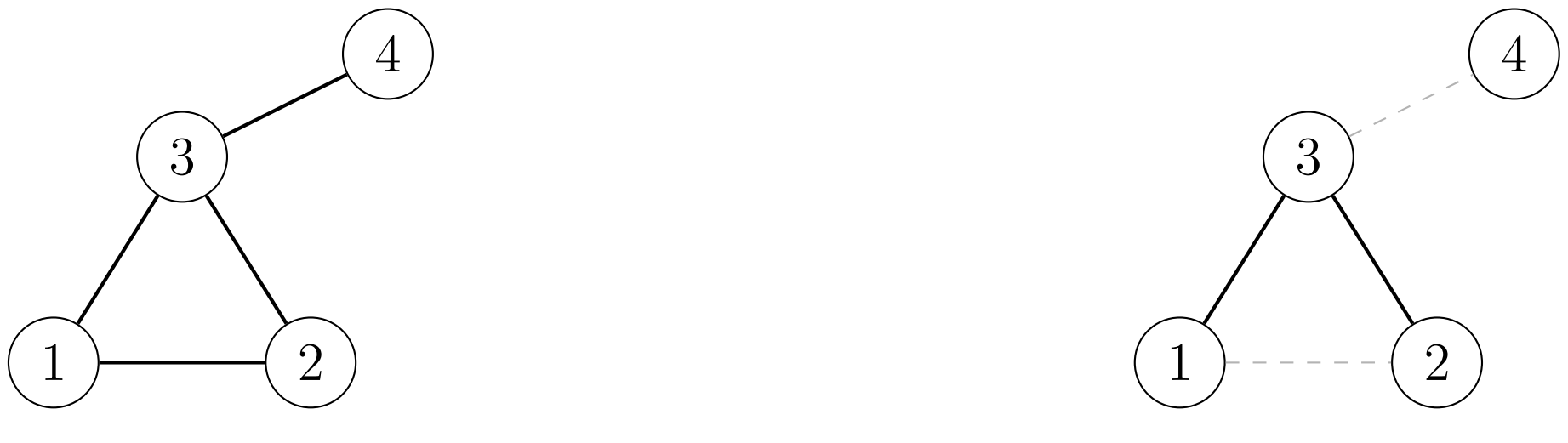
\begin{figure}[H]
\centering
\begin{subfigure}[b]{0.4\textwidth}
\centering
\begin{tikzpicture}[vertex/.style={circle, draw, fill=white, minimum size=7mm, inner sep=0pt}]
\node[vertex] (1) at (0,0)     {$1$};
\node[vertex] (2) at (2,0)     {$2$};
\node[vertex] (3) at (1,1.6)   {$3$};
\node[vertex] (4) at (2.6,2.4) {$4$};
\draw[thick] (1) -- (2);
\draw[thick] (1) -- (3);
\draw[thick] (2) -- (3);
\draw[thick] (3) -- (4);
\end{tikzpicture}
\caption{The graph $G$.}
\end{subfigure}
\hfill
\begin{subfigure}[b]{0.5\textwidth}
\centering
\begin{tikzpicture}[vertex/.style={circle, draw, fill=white, minimum size=7mm, inner sep=0pt}]
\node[vertex] (1) at (0,0)     {$1$};
\node[vertex] (2) at (2,0)     {$2$};
\node[vertex] (3) at (1,1.6)   {$3$};
\node[vertex] (4) at (2.6,2.4) {$4$};
\draw[thick] (1) -- (3);
\draw[thick] (2) -- (3);
\draw[dashed, gray!60] (1) -- (2);
\draw[dashed, gray!60] (3) -- (4);
\end{tikzpicture}
\caption{The $G$-threshold graph, with edge set $\{13,23\}$.}
\end{subfigure}
\caption{A graph $G$ and the $G$-threshold graph obtained from each of the three signed sequences. Dashed edges belong to $G$ but not to the $G$-threshold graph.}
\label{fig:G_and_threshold}
\end{figure}

\medskip
\noindent\textbf{The other two signed sequences give the same graph.}
\begin{itemize}
\item For $(\pi',\epsilon')$, the vertex $1$ is added first, so its sign has no effect. The remaining steps are the same as before, so we get $13$ and $23$, and no other edge.
\item For $(\pi'',\epsilon'')$, we start with $2$. Then $1$ is added with sign $-$, so $12$ is not added. Vertex $3$ has sign $+$, so we add $23$ and $13$. Vertex $4$ has sign $-$, so nothing is added.
\end{itemize}
Thus different signed sequences can give the same $G$-threshold graph.

\end{example}

Observe that every $G$-threshold graph is a subgraph of $G$, but a $G$-threshold graph need not have a unique signed sequence representation (see Example~\ref{ex:threshold_graph}). An edge $ij\in E(G)$ is present in $T_G(\pi,\epsilon)$ if and only if the later endpoint of $ij$ in $\pi$ has sign $+$. So $T_G(\pi,\epsilon)$ depends only on the signs of the later endpoints of the edges of $G$. This motivates the following equivalence relation on signed sequences. A particular case for $G=\mathbb{K}_n$ was introduced by Spiro in \cite{Spiro2019}.

\begin{definition}[$G$-equivalent signed sequences]\label{def:G_equivalent}
Let $G$ be a graph on $[n]$. Two signed sequences $(\pi,\epsilon)$ and $(\sigma,\delta)$ are \textit{$G$-equivalent} if, for every edge $ij\in E(G)$, the later-occurring endpoint of $ij$ has the same sign in $(\pi,\epsilon)$ and in $(\sigma,\delta)$. This is an equivalence relation on $\mathcal S_n$. We denote the set of equivalence classes by $\mathcal S_n(G)$.
\end{definition}

We now show that $G$-equivalence is exactly the relation of giving the same $G$-threshold graph.

\begin{lemma}\label{lem:bij_sig_seq_and_G_threshold_graph}
Two signed sequences give the same $G$-threshold graph if and only if they are $G$-equivalent. Consequently, the map $\mathcal S_n(G)\to\operatorname{Thr}(G)$ sending the class of $(\pi,\epsilon)$ to $T_G(\pi,\epsilon)$ is a bijection.
\end{lemma}

\begin{proof}
By the observation above, an edge $ij\in E(G)$ is present in $T_G(\pi,\epsilon)$ if and only if its later endpoint has sign $+$. So $T_G(\pi,\epsilon)$ and the signs of the later endpoints determine each other, and two signed sequences give the same $G$-threshold graph if and only if they are $G$-equivalent. Hence the map $\mathcal S_n(G)\to\operatorname{Thr}(G)$ is well defined and injective, and it is surjective by the definition of $\operatorname{Thr}(G)$.
\end{proof}

In the rest of this section, we establish a bijection between $\mathcal{S}_n(G)$ and the set of regions of the $G$-threshold arrangement $\mathcal{T}(G)$, which we denote by $\mathcal{R}(\mathcal{T}(G))$. Composing the two bijections then yields the desired correspondence between $\operatorname{Thr}(G)$ and $\mathcal{R}(\mathcal{T}(G))$.

\medskip

\noindent\textbf{The map $\Phi:\mathcal{S}_n\rightarrow \mathcal{R}(\mathcal{T}(G))$:}
For a signed sequence $(\pi,\epsilon)\in\mathcal{S}_n$, define $\Phi(\pi,\epsilon)$ to be the region determined by the following conditions:
\begin{enumerate}
\item If $ij\in E(G)$ and the later-occurring endpoint of $ij$ in $\pi$ has sign $+$, then the region satisfies $x_i+x_j>0$.
\item If $ij\in E(G)$ and the later-occurring endpoint of $ij$ in $\pi$ has sign $-$, then the region satisfies $x_i+x_j<0$.
\end{enumerate}

\medskip

\noindent Before proving that $\Phi$ induces a bijection, we illustrate the map with an example.

\begin{example}\label{ex:Phi}
Let $G$ be the graph from Example~\ref{ex:threshold_graph}, with $E(G)=\{12,\,13,\,23,\,34\}$.
The arrangement $\mathcal{T}(G)$ consists of the hyperplanes
\[
x_1+x_2=0,\quad x_1+x_3=0,\quad x_2+x_3=0,\quad x_3+x_4=0
\]
in $\mathbb{R}^4$. Consider the signed sequence $(\pi,\epsilon)=(1234,\,(+,-,+,-))$.
We determine $\Phi(\pi,\epsilon)$ edge by edge.
\begin{itemize}
\item Edge $12$: the later endpoint is $2$, with sign $-$, so $x_1+x_2<0$.
\item Edge $13$: the later endpoint is $3$, with sign $+$, so $x_1+x_3>0$.
\item Edge $23$: the later endpoint is $3$, with sign $+$, so $x_2+x_3>0$.
\item Edge $34$: the later endpoint is $4$, with sign $-$, so $x_3+x_4<0$.
\end{itemize}
Hence
\[
\Phi(\pi,\epsilon)=\{x\in\mathbb{R}^4 : x_1+x_2<0,\ x_1+x_3>0,\ x_2+x_3>0,\ x_3+x_4<0\}.
\]
This region is nonempty; for instance, $(-1,-1,3,-4)$ lies in it.
Note that the edges $13$ and $23$, which are exactly the edges of $T_G(\pi,\epsilon)$, are the ones with $x_i+x_j>0$.

Now consider $(\pi',\epsilon')=(1234,\,(-,-,+,-))$ and $(\pi'',\epsilon'')=(2134,\,(-,-,+,-))$ from Example~\ref{ex:threshold_graph}.
We will get
\[
\Phi(\pi,\epsilon)=\Phi(\pi',\epsilon')=\Phi(\pi'',\epsilon'').
\]
So the three signed sequences that give the same $G$-threshold graph also give the same region of $\mathcal{T}(G)$.
\end{example}

We first show that the map is well-defined, i.e., that the inequalities defining $\Phi(\pi,\epsilon)$ determine a nonempty region. Given $(\pi,\epsilon)\in\mathcal{S}_n$, define a point $x=(x_1,\ldots,x_n)\in\mathbb{R}^n$ by 
$$
x_{\pi(i)}=\epsilon_i i,
$$ 
for all $ i\in [n]$. Suppose $ij\in E(G)$, and assume that $j=\pi(k)$ occurs later than $i=\pi(\ell)$ in $\pi$, so that $k>\ell$. If $\epsilon_k=+$, then
\[
x_i+x_j\geq -\ell+k>0,
\]
while if $\epsilon_k=-$, then
\[
x_i+x_j\leq \ell-k<0.
\]
Thus $x$ satisfies all the inequalities defining $\Phi(\pi,\epsilon)$. Hence these inequalities have a nonempty intersection, and therefore $\Phi(\pi,\epsilon)$ is well-defined.

\begin{lemma}\label{lem:bij_sig_seq_and_region_G_threshold_arrangement}
The map $\Phi$ maps two signed sequences to the same region if and only if they are $G$-equivalent. Moreover, $\Phi$ is surjective. Therefore, $\Phi$ induces a bijection between the set $\mathcal{S}_n(G)$ and the regions of $\mathcal{T}(G)$.
\end{lemma}

\begin{proof}
Every region of $\mathcal{T}(G)$ lies strictly on one side of each hyperplane in $\mathcal{T}(G)$. By the definition of $\Phi$, for each edge $ij\in E(G)$, the region $\Phi(\pi,\epsilon)$ lies in the half-space $x_i+x_j>0$ or $x_i+x_j<0$ according to whether the later-occurring endpoint of $ij$ in $\pi$ has sign $+$ or $-$. Hence, two signed sequences determine the same region if and only if, for every edge $ij\in E(G)$, the later-occurring endpoints of $ij$ have the same sign in both signed sequences. Equivalently, $(\pi,\epsilon)$ and $(\sigma,\delta)$ are $G$-equivalent.

It remains to prove that $\Phi$ is surjective. Let $\Delta$ be a region of $\mathcal T(G)$ and let $x=(x_1,\dots,x_n)$ be a point in $\Delta$. Let $\pi$ be a permutation that lists the vertices in weakly increasing order of $|x_v|$, with ties broken arbitrarily, and let $\epsilon_k$ be the sign of $x_{\pi_k}$ (if $x_{\pi_k}=0$, choose $\epsilon_k$ arbitrarily). We claim that $\Phi(\pi,\epsilon)=\Delta$.

Let $ij\in E(G)$ with $\pi_r=i$ and $\pi_s=j$, where $r<s$. Then $|x_i|\le|x_j|$. Since $x\in\Delta$, we have $x_i+x_j\neq 0$. This forces $x_j\neq0$, because $x_j=0$ would give $x_i=0$ as well. Therefore $x_i+x_j$ has the same sign as $x_j$, which is $\epsilon_s$. So $x$ satisfies the condition that $(\pi,\epsilon)$ imposes on the hyperplane $x_i+x_j=0$. As this holds for every edge of $G$, the point $x$ lies in $\Phi(\pi,\epsilon)$. Two regions that share a point are equal, so $\Phi(\pi,\epsilon)=\Delta$.

Thus $\Phi$ is surjective. By the first part, its fibers are exactly the $G$-equivalence classes, so $\Phi$ induces a bijection between $\mathcal S_n(G)$ and the regions of $\mathcal T(G)$.
\end{proof}

Combining Lemmas~\ref{lem:bij_sig_seq_and_G_threshold_graph} and~\ref{lem:bij_sig_seq_and_region_G_threshold_arrangement}, we obtain the main result of this section.

\begin{theorem}\label{thm:regions_threshold_graphs}
Let $G$ be a graph on $[n]$. Then the $G$-threshold graphs are in bijection with the regions of $\mathcal T(G)$. In particular, $|\operatorname{Thr}(G)|=|\mathcal R(\mathcal T(G))|$.
\end{theorem}

\begin{proof}
Both $\operatorname{Thr}(G)$ and $\mathcal R(\mathcal T(G))$ are in bijection with $\mathcal S_n(G)$, by Lemma~\ref{lem:bij_sig_seq_and_G_threshold_graph} and Lemma~\ref{lem:bij_sig_seq_and_region_G_threshold_arrangement} respectively. Composing one bijection with the inverse of the other gives the result.
\end{proof}

\section{Faces of the $G$-threshold arrangement}

In this section, we give a description of the faces of the $G$-threshold arrangement. 

\begin{definition}
A \textit{colored signed sequence} is a triple $(\pi,\epsilon,\lambda)$, where $\pi=\pi_1\pi_2\dots\pi_n$ is a permutation of $[n]$, $\epsilon=\epsilon_1\epsilon_2\dots\epsilon_n\in\{+,-,0\}^n$ is a sign sequence, and $\lambda=(\lambda_1,\lambda_2,\dots)$ is a composition of $n$. The composition $\lambda$ assigns color $1$ to the first $\lambda_1$ vertices in $\pi$, color $2$ to the next $\lambda_2$ vertices, and so on. We denote by $f(i)$ the color assigned to vertex $i$. We require the following conditions:
\begin{enumerate}
\item The subgraph induced by each color class is connected.

\item There exists a non-negative integer $k\leq\ell(\lambda)$ such that every vertex in the color classes $1,\ldots,k$ has sign $0$, there are no edges between distinct color classes among these classes, and each of these color classes induces a non-bipartite graph.

\item Each color class with index greater than $k$ contains no vertex with sign $0$. Moreover, if such a color class contains at least two vertices, then every pair of its vertices can be connected by a path in $G$ inside the class along which the signs alternate.
\end{enumerate}
We denote the set of all colored signed sequences by $\mathcal{CS}_n$.
\end{definition}
By conditions 2 and 3, in every color class either all vertices have sign $0$ or none does.

We extend the construction of a $G$-threshold graph by grouping the vertices into color classes. In our convention, a vertex with sign $0$ is regarded as having both signs $+$ and $-$. Thus two vertices of the same color class have \emph{different signs} if one is $+$ and the other is $-$, or if both have sign $0$, and they have the \emph{same sign} $s\in\{+,-\}$ if both have sign $s$.

\begin{definition}[Colored $G$-threshold graph]
Let $(\pi,\epsilon,\lambda)$ be a colored signed sequence, with color function $f$. The \emph{colored $G$-threshold graph} associated to $(\pi,\epsilon,\lambda)$ is the labelled graph on $[n]$, together with its color classes and marked pairs, obtained recursively as follows:
\begin{enumerate}
\item Start with the vertex $\pi_1$.
\item Suppose that $\pi_1,\dots,\pi_{i-1}$ have been added. Add $\pi_i$. For every $k<i$ with $\pi_k\pi_i\in E(G)$:
\begin{itemize}
\item if $f(\pi_k)\neq f(\pi_i)$, join $\pi_k$ and $\pi_i$ if and only if $\pi_i$ has sign $+$;
\item if $f(\pi_k)=f(\pi_i)$, join $\pi_k$ and $\pi_i$ if they have different signs; if they have the same sign $s\in\{+,-\}$, do not join them, but mark the pair $\pi_k\pi_i$ with $s$.
\end{itemize}
\item Color the vertices according to their color classes.
\end{enumerate}
The set of all colored $G$-threshold graphs is denoted by $\mathcal{C}\operatorname{Thr}(G)$.
\end{definition}

Similar to Section~2, we illustrate the construction of a colored $G$-threshold graph with an example. It shows how the graph is built step by step, and that different colored signed sequences can give the same colored $G$-threshold graph.

\begin{example}\label{ex:zero_and_mark}
Let $G$ be the graph on $[7]$ with
\[
E(G)=\{12,13,23,45,46,56,37,67\},
\]
that is, two triangles $\{1,2,3\}$ and $\{4,5,6\}$, and a vertex $7$ adjacent to $3$ and $6$ (see Figure~\ref{fig:zero_and_mark}(a)). Consider the colored signed sequences
\begin{align*}
(\pi,\epsilon,\lambda)&=(1234567,(0,0,0,+,-,-,+),(3,3,1)),\\
(\pi',\epsilon',\lambda')&=(2134567,(0,0,0,+,-,-,+),(3,3,1)),\\
(\pi'',\epsilon'',\lambda'')&=(1234657,(0,0,0,+,-,-,+),(3,3,1)).
\end{align*}
In each of them the color classes are $\{1,2,3\}$, $\{4,5,6\}$ and $\{7\}$. The class $\{1,2,3\}$ is the zero class ($k=1$): all its vertices have sign $0$, and it induces a triangle, which is connected and non-bipartite. The class $\{4,5,6\}$ induces a triangle, so it is connected, it contains no sign $0$, and the vertex $4$ has a sign different from those of $5$ and $6$, so any two of its vertices are joined by a path along which the signs alternate. The class $\{7\}$ consists of one vertex. Hence the conditions of a colored signed sequence hold.

\medskip

\noindent\textbf{Constructing the graph for $(\pi,\epsilon,\lambda)$:}
We start with the vertex $1$. Next, we add $2$ and $3$, both of which have sign $0$. Since the pairs $12,13,23$ are edges of $G$ within the same color class and their endpoints have different signs, we add all three edges. We then add $4$, which has sign $+$. Since it has no earlier neighbor in $G$, no edge is added. Next, we add $5$, which has sign $-$. Its only earlier neighbor in $G$ is $4$, which belongs to the same color class and has a different sign. Thus, we add the edge $45$. We then add $6$, which also has sign $-$. Since $4$ belongs to the same color class and has a different sign, we add the edge $46$. On the other hand, $5$ has the same color and sign as $6$, and $56\in E(G)$, so we mark the pair $56$ with $-$. Finally, we add $7$, which has sign $+$. Its earlier neighbors in $G$ are $3$ and $6$, both of which belong to different color classes from $7$. Thus, we add the edges $37$ and $67$. The resulting graph has edge set $\{12,13,23,45,46,37,67\}$, color classes $\{1,2,3\}$, $\{4,5,6\}$, and $\{7\}$, and the pair $56$ marked with $-$ (see Figure~\ref{fig:zero_and_mark}(b)).

\medskip

\noindent\textbf{The other two sequences give the same graph.}
\begin{itemize}
\item For $(\pi',\epsilon',\lambda')$, the vertices of the zero class are added in the order $2,1,3$. Each of them has sign $0$, and all three pairs are edges of $G$, so we again get $12,13,23$. The other vertices are added exactly as before.
\item For $(\pi'',\epsilon'',\lambda'')$, the vertices $5$ and $6$ are interchanged. They have the same sign, so we again add $45$ and $46$, and the pair $56$ is marked with $-$. Vertex $7$ is added as before.
\end{itemize}
Hence all three sequences give the same colored $G$-threshold graph. What matters is the following. The zero class contributes all the edges of $G$ inside it. In the class $\{4,5,6\}$, the vertex $4$ has a sign different from $5$ and $6$, which gives the edges $45$ and $46$, while $5$ and $6$ have the same sign $-$, which gives the mark on $56$. For the edges $37$ and $67$ between classes, the later endpoint is $7$, whose sign is $+$ in all three sequences, so both are present.

Now consider $(\tau,\theta,\xi)=(1234567,(0,0,0,-,+,+,+),(3,3,1))$. It has the same color classes and the same edges $\{12,13,23,45,46,37,67\}$: the vertex $4$ has a sign different from $5$ and $6$, the later endpoint of $37$ and $67$ is $7$ with sign $+$, and no edge joins $\{4,5,6\}$ to an earlier class. But now $5$ and $6$ both have sign $+$, so the pair $56$ is marked $+$, and the graph is different. A sequence with the sign $-$ on $7$ gives a different set of edges, since $37$ and $67$ are absent.

\begin{figure}[ht]
\centering
\begin{subfigure}[t]{0.47\textwidth}
\centering
\begin{tikzpicture}[scale=0.8,vertex/.style={circle, draw, fill=white, minimum size=7mm, inner sep=0pt}]
\useasboundingbox (-0.6,-0.6) rectangle (6.2,3.2);
\node[vertex] (1) at (0,0)   {$1$};
\node[vertex] (2) at (2,0)   {$2$};
\node[vertex] (3) at (1,1.6) {$3$};
\node[vertex] (4) at (3.2,0) {$4$};
\node[vertex] (5) at (5.2,0) {$5$};
\node[vertex] (6) at (4.2,1.6) {$6$};
\node[vertex] (7) at (2.6,2.8) {$7$};
\draw[thick] (1) -- (2);
\draw[thick] (1) -- (3);
\draw[thick] (2) -- (3);
\draw[thick] (4) -- (5);
\draw[thick] (4) -- (6);
\draw[thick] (5) -- (6);
\draw[thick] (3) -- (7);
\draw[thick] (6) -- (7);
\end{tikzpicture}
\caption{The graph $G$.}
\end{subfigure}
\hfill
\begin{subfigure}[t]{0.47\textwidth}
\centering
\begin{tikzpicture}[scale=0.8,vertex/.style={circle, draw, minimum size=7mm, inner sep=0pt}]
\useasboundingbox (-0.6,-0.6) rectangle (6.2,3.2);
\node[vertex, fill=green!20]  (1) at (0,0)   {$1$};
\node[vertex, fill=green!20]  (2) at (2,0)   {$2$};
\node[vertex, fill=green!20]  (3) at (1,1.6) {$3$};
\node[vertex, fill=blue!15]   (4) at (3.2,0) {$4$};
\node[vertex, fill=blue!15]   (5) at (5.2,0) {$5$};
\node[vertex, fill=blue!15]   (6) at (4.2,1.6) {$6$};
\node[vertex, fill=orange!25] (7) at (2.6,2.8) {$7$};
\draw[thick] (1) -- (2);
\draw[thick] (1) -- (3);
\draw[thick] (2) -- (3);
\draw[thick] (4) -- (5);
\draw[thick] (4) -- (6);
\draw[dotted, thick] (5) -- (6);
\node at (4.95,0.95) {$-$};
\draw[thick] (3) -- (7);
\draw[thick] (6) -- (7);
\end{tikzpicture}
\caption{The colored $G$-threshold graph.}
\end{subfigure}
\caption{A graph $G$ and the colored $G$-threshold graph of Example~\ref{ex:zero_and_mark}, with the zero class $\{1,2,3\}$ (green), the class $\{4,5,6\}$ (blue) and the class $\{7\}$ (orange). The dotted pair is marked $-$.}
\label{fig:zero_and_mark}
\end{figure}
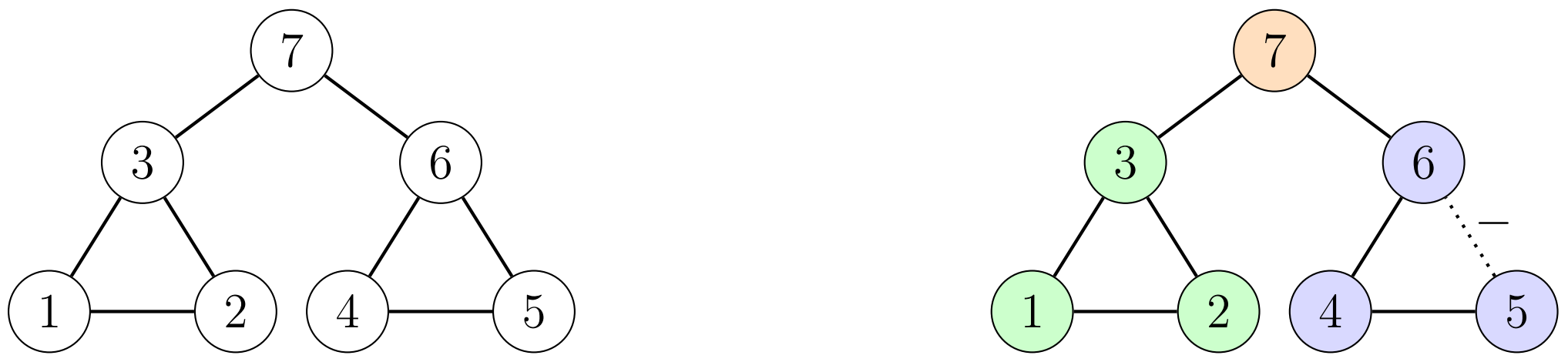
\end{example}
Again, distinct colored signed sequence may give rise to the same colored $G$-threshold graph. To tackle this issue, we define the following equivalence relation on colored signed sequences.

\begin{definition}[$G$-equivalent colored signed sequences]\label{def:G_equivalent_colored}
Let $G$ be a graph on $[n]$. For a colored signed sequence and an edge $ij\in E(G)$ whose endpoints lie in the same color class, we say that $ij$ has \emph{type} $0$ if $i$ and $j$ have different signs, type $+$ if both have sign $+$, and type $-$ if both have sign $-$ (a pair of vertices of sign $0$ has different signs). Two colored signed sequences $(\pi,\epsilon,\lambda_1)$ and $(\sigma,\delta,\lambda_2)$ are \emph{$G$-equivalent} if the following hold.
\begin{enumerate}
\item The two sequences have the same color classes, that is, $i$ and $j$ lie in the same color class in $\lambda_1$ if and only if they lie in the same color class in $\lambda_2$, for all $i,j\in[n]$.
\item Every edge $ij\in E(G)$ inside a color class has the same type in both sequences.
\item For every edge $ij\in E(G)$ joining different color classes, the later-occurring endpoint of $ij$ has the same sign in $(\pi,\epsilon,\lambda_1)$ and in $(\sigma,\delta,\lambda_2)$.
\end{enumerate}
This is an equivalence relation on $\mathcal{CS}_n$. We denote the set of equivalence classes by $\mathcal{CS}_n(G)$.
\end{definition}

We now extend Lemma~\ref{lem:bij_sig_seq_and_G_threshold_graph} to colored signed sequences and colored $G$-threshold graphs.

\begin{lemma}\label{lem:bij_colored_sig_seq_and_G_threshold_graph}
Two colored signed sequences give the same colored $G$-threshold graph if and only if they are $G$-equivalent. Consequently, the map $\mathcal{CS}_n(G)\to\mathcal{C}\mathrm{Thr}(G)$ sending the $G$-equivalence class of $(\pi,\epsilon,\lambda)$ to $T_G(\pi,\epsilon,\lambda)$ is a bijection.
\end{lemma}

\begin{proof}
The proof is similar to that of Lemma~\ref{lem:bij_sig_seq_and_G_threshold_graph}. Let $T$ be the colored $G$-threshold graph of $(\pi,\epsilon,\lambda)$. By the definition of $T$, it records the following for every edge $ij\in E(G)$.
\begin{itemize}
\item If $i$ and $j$ belong to the same color class, then $ij\in E(T)$ if and only if $i$ and $j$ have different signs (a pair of vertices of sign $0$ has different signs), and otherwise $ij$ is marked $+$ if $i$ and $j$ both have sign $+$, and marked $-$ if both have sign $-$.
\item If $i$ and $j$ belong to different color classes, then $ij\in E(T)$ if and only if the later endpoint of $ij$ has sign $+$.
\end{itemize}
Moreover, $T$ determines the partition of $[n]$ into color classes. Hence two colored signed sequences give the same colored $G$-threshold graph if and only if the following hold: their partitions into color classes are the same; for every edge $ij\in E(G)$ inside a color class, $i$ and $j$ are both positive, both negative, or have different signs in one sequence if and only if the same holds in the other; and for every edge between different color classes, the later endpoints have the same sign in both. These are precisely the conditions for the two colored signed sequences to be $G$-equivalent.

Consequently, the map sending the $G$-equivalence class of $(\pi,\epsilon,\lambda)$ to $T_G(\pi,\epsilon,\lambda)$ is well defined and injective. It is surjective by the definition of $\mathcal{C}\mathrm{Thr}(G)$, and hence is a bijection.
\end{proof}

In the rest of this section, we establish a bijection between $\mathcal{CS}_n(G)$ and the set of faces of the $G$-threshold arrangement $\mathcal{T}(G)$, which we denote by $\mathcal{F}(\mathcal{T}(G))$. Composing the two bijections then yields the desired correspondence between $\mathcal{C}\operatorname{Thr}(G)$ and $\mathcal{F}(\mathcal{T}(G))$.

\medskip

\noindent\textbf{The map $\tilde\Phi:\mathcal{CS}_n\rightarrow\mathcal{F}(\mathcal{T}(G))$.}
For a colored signed sequence $(\pi,\epsilon,\lambda)\in\mathcal{CS}_n$, define $\tilde\Phi(\pi,\epsilon,\lambda)$ to be the set of points $\mathbf{x}\in\mathbb{R}^n$ satisfying the following conditions for every edge $ij\in E(G)$.
\begin{enumerate}
\item If $i$ and $j$ lie in different color classes and the later-occurring endpoint of $ij$ in $\pi$ has sign $+$, then $x_i+x_j>0$.
\item If $i$ and $j$ lie in different color classes and the later-occurring endpoint of $ij$ in $\pi$ has sign $-$, then $x_i+x_j<0$.
\item If $i$ and $j$ lie in the same color class, then
\begin{itemize}
\item $x_i+x_j=0$ if $ij$ has type $0$,
\item $x_i+x_j>0$ if $ij$ has type $+$,
\item $x_i+x_j<0$ if $ij$ has type $-$.
\end{itemize}
\end{enumerate}

\medskip

\noindent Before proving that $\tilde\Phi$ induces a bijection, we illustrate the map with an example.

\begin{example}
Let $G$ be the graph from Example~\ref{ex:zero_and_mark}. We determine $\tilde\Phi$ for the colored signed sequence $(\pi,\epsilon,\lambda) =(1234567, (0,0,0,+,-,-,+),(3,3,1))$.

\begin{itemize}
\item The edges $12,13,23$ lie inside the zero class and have type $0$, so $x_1=x_2=x_3=0$.
\item The edges $45,46$ have type $0$, so $x_5=x_6=-x_4$. Write $x_4=t$. The edge $56$ has type $-$, so $x_5+x_6<0$, that is, $t>0$.
\item The edges $37$ and $67$ join different classes, and their later endpoint $7$ has sign $+$, so $x_3+x_7>0$ and $x_6+x_7>0$. Writing $x_7=u$, this gives $u>0$ and $u>t$.
\end{itemize}
Hence
\[
\tilde\Phi(\pi,\epsilon,\lambda)=\{\mathbf x\in\mathbb R^7\mid\ x_1=x_2=x_3=0,\ x_4=-x_5=-x_6=t,\ x_7=u,\ 0<t<u\}.
\]
For instance, $(0,0,0,1,-1,-1,2)$ lies in it. Its dimension is $2$, in agreement with Theorem~\ref{thm:dim_through_thresh_grph}: the zero class is non-bipartite and contributes nothing, and the classes $\{4,5,6\}$ and $\{7\}$ are bipartite and contribute one parameter each. The sequence $(\tau,\theta,\xi)$ from Example~\ref{ex:zero_and_mark} gives the face
\[
\{\mathbf x\in\mathbb R^7:\ (x_1,x_2,x_3)=(0,0,0),\ (x_4,x_5,x_6)=(t,-t,-t),\ x_7=u,\ t<0,\ u>0\},
\]
since now $x_5+x_6>0$ gives $t<0$, and $x_3+x_7>0$ gives $u>0$, while $x_6+x_7>0$ holds automatically. It contains, for instance, $(0,0,0,-1,1,1,2)$. The two faces have the same solid edges in the colored graph, and they differ only in the mark on $56$, which separates $t>0$ from $t<0$. This shows that the marks are needed.
\end{example}

\noindent We first show that the map $\tilde\Phi$ is well-defined, i.e., that the conditions defining $\tilde\Phi(\pi,\epsilon,\lambda)$ determine a nonempty set. Given a colored signed sequence $(\pi,\epsilon,\lambda)\in\mathcal{CS}_n$, define a point $\mathbf{x}=(x_1,\ldots,x_n)\in\mathbb{R}^n$ by
\[
x_{\pi_i}=\epsilon_i\,j,\qquad\text{for all } i\in[n],
\]
where $j$ is the smallest index such that $\pi_j$ belongs to the same color class as $\pi_i$. Thus all vertices of a color class have the same absolute value $j$, or all have the value $0$ if the class is a zero class. Suppose $ij\in E(G)$.

If $i$ and $j$ belong to different color classes, assume that $j$ occurs later than $i$ in $\pi$. Then $\epsilon_j\neq0$: otherwise $j$ would lie in a zero class, and so would $i$, since the zero classes come first, but there are no edges between distinct zero classes. Since the color classes are consecutive blocks in $\pi$, the first vertex of the class of $i$ occurs before the first vertex of the class of $j$. Hence $|x_i|<|x_j|$, and therefore
\[
\operatorname{sgn}(x_i+x_j)=\operatorname{sgn}(x_j)=\epsilon_j .
\]
Thus conditions (1) and (2) are satisfied.

If $i$ and $j$ belong to the same color class, then $x_i$ and $x_j$ have the same absolute value. If the class is a zero class, then $x_i=x_j=0$, and $ij$ has type $0$, as required. Otherwise, no vertex of the class has sign $0$, and $x_i+x_j=0$ if $i$ and $j$ have different signs, $x_i+x_j>0$ if both have sign $+$, and $x_i+x_j<0$ if both have sign $-$. Thus condition (3) is satisfied.

Hence $\mathbf{x}\in\tilde\Phi(\pi,\epsilon,\lambda)$. So the conditions have a nonempty solution set, and they prescribe the sign of $x_i+x_j$ for every edge $ij\in E(G)$, so $\tilde\Phi(\pi,\epsilon,\lambda)$ is a face of $\mathcal{T}(G)$ and $\tilde\Phi$ is well-defined.

\begin{lemma}\label{lem:bij_colored_sig_seq_and_faces_G_threshold_arrangement}
The map $\tilde\Phi$ maps two colored signed sequences to the same face if and only if they are $G$-equivalent. Moreover, $\tilde\Phi$ is surjective. Therefore, $\tilde\Phi$ induces a bijection between the set $\mathcal{CS}_n(G)$ and the faces of $\mathcal{T}(G)$.
\end{lemma}

\begin{proof}
A face of $\mathcal T(G)$ is determined by its sign vector, that is, by the sign of $x_i+x_j$ for every edge $ij\in E(G)$. By the definition of $\tilde\Phi$, the face $\tilde\Phi(\pi,\epsilon,\lambda)$ has the following signs.
\begin{itemize}
\item If $i$ and $j$ belong to different color classes, the sign of $x_i+x_j$ is the sign of the later-occurring endpoint of $ij$.
\item If $i$ and $j$ belong to the same color class, then $x_i+x_j=0$ if $i$ and $j$ have different signs, $x_i+x_j>0$ if both have sign $+$, and $x_i+x_j<0$ if both have sign $-$.
\end{itemize}
Moreover, the face determines the partition of $[n]$ into color classes. Indeed, these classes are precisely the connected components of the graph on $[n]$ with edge set $\{ij\in E(G):x_i+x_j=0\text{ on the face}\}$. To see this, observe that conditions 1--3 of a colored signed sequence ensure that every color class is connected in this graph: in a nonzero class, connectivity follows from alternating paths, whose edges join vertices of different signs and hence have type $0$, while in a zero class, all edges of the induced subgraph have type $0$ and the class is connected. Moreover, no edge of this graph joins vertices from different color classes, since for such an edge $x_i+x_j$ is strictly positive or strictly negative by conditions (1) and (2). Hence two colored signed sequences determine the same face if and only if their partitions into color classes agree, every edge inside a color class has the same type in both sequences, and for every edge between different color classes the later endpoints have the same sign in both. These are precisely the conditions for the two sequences to be $G$-equivalent.

It remains to prove that $\tilde\Phi$ is surjective. Let $\Delta$ be a face of $\mathcal T(G)$. Define a relation $\sim$ on $[n]$ by declaring $i\sim j$ if there is a path from $i$ to $j$ in $G$ all of whose edges $ab$ satisfy $x_a+x_b=0$ on $\Delta$. This is an equivalence relation, and the equivalence classes are the color classes. Along such a path the coordinates are $\pm$ each other, so all vertices in a color class have the same absolute value, and a class containing one zero coordinate consists of zero coordinates only. Choose $\mathbf x=(x_1,\ldots,x_n)\in\Delta$ with the minimum possible number of zero coordinates.

We claim that these classes satisfy the conditions in the definition of a colored signed sequence. Each class is connected in $G$ by the paths defining $\sim$. Consider a class $C$ with all coordinates zero. Every edge of $G$ inside $C$ then satisfies $x_i+x_j=0$, so no two distinct such classes are joined by an edge of $G$ (the endpoints would be $\sim$-equivalent). Suppose that $C$ induces a bipartite graph with parts $A\sqcup B$. Replace the coordinates by $\delta$ on $A$ and $-\delta$ on $B$, with $\delta>0$ small. All the equalities inside $C$ are preserved, and for an edge from $C$ to a vertex $j$ outside $C$ we have $x_j\neq0$ by the previous remark, so the sign of $x_i+x_j$ is preserved for small $\delta$. The new point lies in $\Delta$ and has fewer zero coordinates, contrary to the choice of $\mathbf x$. So $C$ induces a non-bipartite graph. A class containing no zero coordinate has nonzero coordinates, and along each equality edge inside it the signs are opposite, so any two of its vertices are joined by a path along which the signs alternate.

Order the color classes by increasing absolute value of their coordinates, breaking ties arbitrarily; the zero classes come first. Order the vertices within each class arbitrarily. This gives a permutation $\pi$ in which the vertices of each color class occur consecutively. Let $\lambda$ record the sizes of the color classes in this order, and let $\epsilon_k$ be the sign of $x_{\pi_k}$, with $\epsilon_k=0$ when $x_{\pi_k}=0$. Then $(\pi,\epsilon,\lambda)$ is a colored signed sequence.

Now consider an edge $ij\in E(G)$. If $i$ and $j$ belong to the same color class, then $x_i+x_j=0$ if their signs differ or both are $0$, $x_i+x_j>0$ if both are positive, and $x_i+x_j<0$ if both are negative, as required. Otherwise, suppose that $i$ and $j$ belong to different color classes, with $i$ occurring before $j$ in $\pi$. By the ordering of the color classes, $|x_i|\leq|x_j|$. Since $i$ and $j$ belong to different color classes, $x_i+x_j\neq0$ on $\Delta$. It follows that $x_j\neq0$ and that $x_i+x_j$ has the same sign as $x_j$, namely $\epsilon_j$. Hence $\mathbf{x}$ satisfies all the conditions defining $\tilde\Phi(\pi,\epsilon,\lambda)$, so $\mathbf{x}$ lies in this face, and therefore
\[
\tilde\Phi(\pi,\epsilon,\lambda)=\Delta.
\]
Thus $\tilde\Phi$ is surjective. Since its fibers are exactly the $G$-equivalence classes, it induces a bijection between $\mathcal{CS}_n(G)$ and the faces of $\mathcal T(G)$.
\end{proof}

Combining Lemmas~\ref{lem:bij_colored_sig_seq_and_G_threshold_graph} and~\ref{lem:bij_colored_sig_seq_and_faces_G_threshold_arrangement}, we obtain the main result of this section.

\begin{theorem}\label{thm:faces_threshold_graphs}
Let $G$ be a graph on $[n]$. Then the colored $G$-threshold graphs are in bijection with the faces of $\mathcal T(G)$. In particular, $|\mathcal{C}\operatorname{Thr}(G)|=|\mathcal F(\mathcal T(G))|$.
\end{theorem}

\begin{proof}
Both $\mathcal{C}\operatorname{Thr}(G)$ and $\mathcal F(\mathcal T(G))$ are in bijection with $\mathcal{C}\mathcal S_n(G)$, by Lemma~\ref{lem:bij_colored_sig_seq_and_G_threshold_graph} and Lemma~\ref{lem:bij_colored_sig_seq_and_faces_G_threshold_arrangement} respectively. Composing one bijection with the inverse of the other gives the result.
\end{proof}

\subsection*{Dimension of faces of $G$-threshold arrangement}

In the next theorem, we show how the dimension of a face of the $G$-threshold arrangement can be read off from its colored $G$-threshold graph.

\begin{theorem}\label{thm:dim_through_thresh_grph}
Let $\Delta$ be the face corresponding to a colored $G$-threshold graph $G_0$. Then $\dim(\Delta)$ is equal to the number of bipartite connected components of the graph obtained from $G_0$ by keeping only the edges whose endpoints lie in the same color class.
\end{theorem}

\begin{proof}
Let $H$ denote the graph obtained from $G_0$ by keeping only the edges whose endpoints lie in the same color class (marked pairs are not edges). Each color class is connected in $H$: in a nonzero class, any two vertices are joined by a path in $G$ along which the signs alternate, and every edge of such a path joins vertices of different signs, hence lies in $H$; in a zero class, every edge of $G$ inside it joins two vertices of sign $0$, hence lies in $H$, and the class is connected in $G$ by condition~1. No edge of $H$ joins two different color classes. Hence the connected components of $H$ are exactly the color classes.

The equality relations defining $\Delta$ arise from the edges of $H$: an edge $ij\in E(G)$ with $i,j$ in the same color class satisfies $x_i+x_j=0$ on $\Delta$ if and only if $ij\in E(H)$. All other conditions defining $\Delta$ are strict inequalities. Let
\[
L=\{x\in\mathbb R^n : x_i+x_j=0 \text{ for every edge } ij \text{ of } H\}.
\]
Then $\Delta$ is a nonempty subset of $L$ defined by finitely many strict linear inequalities, and hence it is open in $L$, and hence $\dim(\Delta)=\dim(L)$. We show that each bipartite component of $H$ contributes one free parameter to $L$, whereas each non-bipartite component contributes none.

First, consider a non-bipartite component $C$ of $H$, that is, one containing an odd cycle. Let $i_1,i_2,\ldots,i_{2r+1}$ be the vertices of such a cycle, with indices taken cyclically. The equality relations along the edges of the cycle give
\[
x_{i_1}=-x_{i_2}=x_{i_3}=\cdots=x_{i_{2r+1}}=-x_{i_1}.
\]
Hence $x_{i_1}=0$, and consequently every coordinate corresponding to a vertex on the cycle is zero. Since $C$ is connected, every other vertex of $C$ is joined to the cycle by a path of equalities $x_a+x_b=0$, so its coordinate is also zero. Thus all coordinates in $C$ are zero, and $C$ contributes no free parameter.

Now consider a bipartite component $C$ of $H$ with bipartition $A\sqcup B$. The equality relations along the edges of $C$ say that adjacent vertices have opposite coordinates. Since $C$ is connected, all coordinates corresponding to vertices in $A$ are equal to some value $t$, and all coordinates corresponding to vertices in $B$ are equal to $-t$. Conversely, every value of $t$ satisfies these relations. Thus $C$ contributes exactly one free parameter.

Different components of $H$ involve disjoint sets of coordinates, so they share no equality relations. The edges of $G$ between different color classes only impose strict inequalities, which do not lower the dimension of $\Delta$. Therefore, $\dim(L)$ is the number of bipartite components of $H$, and hence so is $\dim(\Delta)$, as required.
\end{proof}

\section{Counting faces of $\mathcal{T}(\mathbb{K}_n)$ and $\mathcal{T}(\mathbb{K}_{m,n})$}

By Section~3, the $d$-dimensional faces of $\mathcal{T}(G)$ correspond to colored $G$-threshold graphs whose color classes give dimension $d$ (Theorem~\ref{thm:dim_through_thresh_grph}). In this section we count these graphs when $G=\mathbb{K}_n$ and when $G=\mathbb{K}_{m,n}$, by first choosing the color classes and then counting the ways to place the edges between them.

\subsection{When $G=\mathbb{K}_n$}
A \emph{strong composition} of $k$ is a sequence $\lambda=(\lambda_1,\lambda_2,\ldots,\lambda_d)$ of positive integers with $\lambda_1+\cdots+\lambda_d=k$, as opposed to a weak composition, which allows zero parts \cite[Section~1.2]{Stanley2011}. The \emph{Fubini number} $a(r)$ is the number of ordered set partitions of $[r]$, that is,
\[
a(r)=\sum_{j=0}^{r}j!\,S(r,j),
\]
where $S(r,j)$ denotes the Stirling numbers of the second kind; so $a(0),a(1),a(2),\dots=1,1,3,13,75,\dots$ \cite[Section~1.4]{Stanley2011}. Recall that the number of labelled threshold graphs on $[r]$ is
\[
|\mathrm{Thr}(\mathbb K_r)|=2\bigl(a(r)-r\,a(r-1)\bigr)\quad (r\ge 2),\qquad |\mathrm{Thr}(\mathbb K_1)|=1,
\]
which gives $2,8,46,\dots$ for $r=2,3,4$; see \cite{Spiro2019,DeshpandeMenonSingh2021}.

For a strong composition $\lambda=(\lambda_1,\lambda_2,\ldots,\lambda_d)$, let $\mathcal R(\lambda)$ denote the multiset of maximal runs of parts equal to $1$. For $R\in\mathcal R(\lambda)$, let $|R|$ denote the length of the run. We define $\widehat{\lambda}=(\widehat\lambda_1,\ldots,\widehat\lambda_t)$ to be the composition obtained from $\lambda$ by replacing each maximal run of $1$'s by its length and retaining every part of size at least $2$.

For example, if $\lambda=(3,1,1,2,1,1,1,4,1,1)$, then $\mathcal R(\lambda)$ consists of the runs $(1,1)$, $(1,1,1)$, and $(1,1)$. The corresponding composition $\widehat{\lambda}$ is obtained by replacing each maximal run of $1$'s by its length; thus,
\[
\widehat{\lambda}=(3,2,2,3,4,2).
\]

For a strong composition $\lambda$ of $k$, define
$$
C(\lambda,k)=
\begin{cases}
|\mathrm{Thr}(\mathbb K_r)| \displaystyle\prod_{\substack{R\in\mathcal R(\lambda)\\R\neq R_1}} 2a(|R|), & \text{if } k=n \text{ and } \lambda_1=1,\\[1.2em]
\displaystyle\prod_{R\in\mathcal R(\lambda)}2a(|R|), & \text{otherwise},
\end{cases}
$$
where, in the first case, $R_1$ denotes the initial run of $1$'s and $r=|R_1|$. Also define
$$
B(\lambda,k)=
\begin{cases}
\displaystyle\prod_{\substack{2\le i\le d\\\lambda_i\ge2}}(2^{\lambda_i}-2), & \text{if } k=n \text{ and } \lambda_1=2,\\[1.2em]
\displaystyle\prod_{\substack{1\le i\le d\\\lambda_i\ge2}}(2^{\lambda_i}-2), & \text{otherwise}.
\end{cases}
$$
The exceptional cases occur only for $k=n$, that is, when there is no zero class: if $k<n$, the zero class precedes all the other color classes, and since every vertex of $\mathbb K_n$ is adjacent to all the others, the edges to the zero class record the sign of every vertex of the other classes.

\begin{theorem}\label{thm:faces_of_thr_arrang}
For $d\in[n]$, the number of $d$-dimensional faces of the threshold arrangement $\mathcal T(\mathbb K_n)$ is
$$
\sum_{\substack{k\in[n-3]\cup\{n\}\\ \lambda_1+\cdots+\lambda_d=k}}
\binom{n}{\widehat\lambda_1,\ldots,\widehat\lambda_t,\,n-k}\,
C(\lambda,k)\,B(\lambda,k),
$$
where, for each $k$, the inner sum runs over all strong compositions $\lambda=(\lambda_1,\ldots,\lambda_d)$ of $k$, and $t$ denotes the number of parts of $\widehat\lambda$.
\end{theorem}

\begin{proof}
By Theorem~\ref{thm:dim_through_thresh_grph}, to count the $d$-dimensional faces of $\mathcal T(\mathbb K_n)$ we count the colored $\mathbb K_n$-threshold graphs with exactly $d$ bipartite color classes, that is, with $d$ color classes whose vertices have nonzero sign.

Observe that for $\mathbb K_n$ there can be at most one color class containing vertices of sign $0$, since two distinct such color classes would have an edge between them, contradicting the requirement that the zero-sign color classes have no edges between them. Thus the vertices of sign $0$ form a single color class, whose induced graph is the complete graph $\mathbb K_{n-k}$. This class must be non-bipartite, so it contains an odd cycle, and hence $n-k=0$ or $n-k\ge 3$. Equivalently, $k\in[n-3]\cup\{n\}$, where $k=n$ means that there is no zero class.

The remaining $k$ vertices are divided into the $d$ bipartite color classes, whose sizes are recorded by a strong composition $\lambda=(\lambda_1,\ldots,\lambda_d)$ of $k$. We first choose the vertex sets. Each maximal run of $1$'s is regarded as a single block, and its internal structure is counted later. The vertex sets of the blocks of $\widehat\lambda$, together with the zero-sign class, can be chosen in
\[
\binom{n}{\widehat\lambda_1,\ldots,\widehat\lambda_t,\,n-k}
\]
ways.

It remains to count the possible colored $\mathbb K_n$-threshold graphs for a fixed $\lambda$ and a fixed choice of vertex sets. The relative order of the color classes of size at least two is fixed, while each run of $1$'s sits at a fixed position relative to these classes. Indeed, let $B$ and $C$ be color classes with $C$ of size at least two. The edges between different color classes are determined by the signs of the vertices in the later class, and every nonzero color class of size at least two contains vertices of both signs. So the graph tells which of $B$ and $C$ is later. Hence we count the possibilities for each run and for each class of size at least two separately.

Consider a run of $1$'s of length $r$. Each part of the run is a singleton color class, so the vertices of the run can be ordered in different ways, and their signs determine the edges among these vertices. These choices are the standard-form descriptions of threshold graphs on $r$ vertices.

If the run is the initial run, which can occur only when $k=n$ and $\lambda_1=1$, there is no preceding color class. In this case the choices are exactly the labelled threshold graphs on $r$ vertices, giving
$$
|\mathrm{Thr}(\mathbb K_r)|.
$$
If the run occurs after another color class, the first singleton vertex of the run interacts with the preceding color class. In addition to the threshold-pair data describing the run, there are two choices according to the sign of this first vertex. Indeed, since there is an earlier class, all the signs of the vertices of the run are recorded by the edges to it. Two vertices of the run with the same sign can be interchanged without changing the graph. So the run is determined by an ordered set partition of its vertices into blocks of constant sign, with alternating signs, together with the sign of the first block. Thus a run of length $r$ contributes
\[
2a(r)
\]
possibilities. Multiplying these contributions over all maximal runs of $1$'s gives precisely the factor $C(\lambda,k)$.

Next consider a color class of size $\lambda_i\ge2$. Since the class is bipartite, its two parts are determined by the $+$ and $-$ signs. We may choose any nonconstant sign assignment, giving
\[
2^{\lambda_i}-2
\]
possibilities, since the two constant assignments do not give a bipartition. This count is valid whenever the class occurs after an earlier color class, that is, if $k<n$ or $i\ge2$, because then all signs are recorded by the edges to the earlier classes.

There is one exceptional case: $k=n$ and $i=1$, where the first class is the first class of the whole sequence. If $\lambda_1=2$, the two nonconstant sign assignments $(+,-)$ and $(-,+)$ give the same colored threshold graph, so there is exactly one possibility, namely $2^{\lambda_1-1}-1=1$. If $\lambda_1\ge3$, any nonconstant assignment has two vertices of the same sign in the class. This pair is marked with this sign, and exchanging $+$ and $-$ changes the mark, which gives a different face. Hence all $2^{\lambda_1}-2$ assignments are distinct.

Therefore the contribution from the color classes of size at least two is $B(\lambda,k)$. Multiplying the choices for the vertex sets, the runs of $1$'s, and the bipartite color classes of size at least two, we obtain
$$
\binom{n}{\widehat\lambda_1,\ldots,\widehat\lambda_t,\,n-k}\,C(\lambda,k)\,B(\lambda,k)
$$
colored $\mathbb K_n$-threshold graphs for the fixed pair $(k,\lambda)$. Summing over all admissible $k$ and all strong compositions $\lambda$ of $k$ with $d$ parts gives the stated formula.
\end{proof}

\subsection{When $G=\mathbb K_{m,n}$}

Let $G=\mathbb K_{m,n}$ with parts $X$ and $Y$, where $|X|=m$ and $|Y|=n$. 
For a colored signed sequence $(\pi,\epsilon,\lambda)$, the color classes are consecutive blocks of $\pi$, since $\lambda$ is a composition. Hence for two different color classes $B$ and $C$, either all vertices of $B$ occur before all vertices of $C$ in $\pi$, or the other way round. In the first case we say that $B$ is \emph{earlier} than $C$. 
For a set $\mathcal C$ of color classes, let $H_{\mathcal C}$ be the graph whose vertices are the classes in $\mathcal C$, where two classes are adjacent if some edge of $G$ joins them.

\begin{lemma}\label{lem:orient}
Let $\mathcal C$ be a set of color classes of a colored signed sequence for $\mathbb{K}_{m,n}$, with $a$ classes consisting of one vertex of $X$, $b$ classes consisting of one vertex of $Y$, and $p$ classes with at least two vertices. Then $H_{\mathcal C}$ is the join of the complete graph $\mathbb{K}_p$ and $\mathbb{K}_{a,b}$, and the colored $\mathbb K_{m,n}$-threshold graphs with set of color classes $\mathcal C$ are in bijection with the acyclic orientations of $H_{\mathcal C}$.
\end{lemma}

\begin{proof}
Since $\mathbb K_{m,n}$ is bipartite, no color class can induce a non-bipartite graph. Hence no vertex has sign $0$. A color class with at least two vertices is connected, so it contains vertices of both $X$ and $Y$. Any two vertices of the class are joined by a path along which the signs alternate; a path between two vertices of the same part has even length, and one between vertices of different parts has odd length. Any two vertices of the class are joined by a path along which the signs alternate; a path between two vertices of the same part has even length, and one between vertices of different parts has odd length. Hence all vertices of the class in $X$ have a common sign $t$, and all vertices of the class in $Y$ have the sign $-t$. Thus every pair $uw$ with $u\in X$, $w\in Y$ inside the class has different signs, so $uw$ is an edge of the colored $\mathbb K_{m,n}$-threshold graph, and no pair inside a class is marked.

Define the adjusted sign of a vertex to be its sign if it is in $X$, and the opposite of its sign if it is in $Y$. By the previous paragraph, all vertices of a class $C$ with at least two vertices have the same adjusted sign, and this holds trivially for a class with one vertex. We denote it by $t_C$. Two classes with at least two vertices are adjacent in $H_{\mathcal C}$, since each of them contains vertices of both $X$ and $Y$. Such a class is adjacent to every class with one vertex. A class $\{u\}$ with $u\in X$ is adjacent to a class $\{w\}$ with $w\in Y$, and two classes with one vertex in the same part are not adjacent. Thus $H_{\mathcal C}=\mathbb{K}_p\vee\mathbb{K}_{a,b}$.

We now define the map from the colored $\mathbb K_{m,n}$-threshold graphs with set of color classes $\mathcal C$ to the orientations of $H_{\mathcal C}$. Let $T$ be such a graph, let $B$ and $C$ be adjacent classes with $B$ earlier than $C$, and write $X_D=D\cap X$, $Y_D=D\cap Y$ for a class $D$. We say that $C$ wins over $B$ if the edges of $G$ joining $X_C$ and $Y_B$ are in $T$ and those joining $X_B$ and $Y_C$ are not, and that $C$ loses to $B$ in the opposite case, and we orient the edge of $H_{\mathcal C}$ between $B$ and $C$ towards the winner. This is well defined and exactly one case occurs: every edge between $B$ and $C$ is in $T$ if and only if its endpoint in $C$ has sign $+$, and this endpoint has sign $t_C$ if it lies in $X_C$ and $-t_C$ if it lies in $Y_C$. So $C$ wins if $t_C=+$ and loses if $t_C=-$. Hence for every class $C$, all arrows between $C$ and the classes earlier than $C$ point into $C$, or all of them point out of $C$. The orientation is acyclic: in a directed cycle, the class occurring last in $\pi$ would have one incoming and one outgoing arrow, both to earlier classes, which is impossible.

The map is injective, since $T$ is determined by $\mathcal C$ and $O$: the edges of $T$ inside the classes are fixed, and for $u\in X$, $w\in Y$ in different classes, $uw\in E(T)$ if and only if $O$ orients the edge between their classes from the class of $w$ to the class of $u$. For surjectivity, let $O$ be an acyclic orientation of $H_{\mathcal C}$, and list the classes as $B_1,\dots,B_d$ so that every arrow of $O$ goes from $B_i$ to $B_j$ with $i<j$. Let $\pi$ list the vertices of $B_1$, then those of $B_2$, and so on, let $\lambda=(|B_1|,\dots,|B_d|)$, and assign the sign $+$ to the vertices of $X$ and the sign $-$ to the vertices of $Y$. Then $(\pi,\epsilon,\lambda)$ is a colored signed sequence with set of color classes $\mathcal C$. Every class $C$ has $t_C=+$, so $C$ wins over every earlier adjacent class, and the arrow between them goes from the earlier class to the later one. Hence the orientation of the colored graph of this sequence is $O$. This gives the required bijection.
\end{proof}

By Lemma~\ref{lem:orient}, counting the faces of $\mathcal T(\mathbb K_{m,n})$ of dimension $d$ reduces to counting acyclic orientations of the graphs $H_{\mathcal C}$, which we do with chromatic polynomials. Let $S(a,i)$ denote the Stirling numbers of the second kind, $s^{\overline{i}}=s(s+1)\cdots(s+i-1)$, and
$$
c_s(a,b)=\sum_{i=0}^{a}(-1)^{a-i}\,S(a,i)\,s^{\overline{i}}\,(s+i)^{b}.
$$

\begin{theorem}\label{thm:faces_Kmn}
For $d\ge 1$, the number of $d$-dimensional faces of the threshold arrangement $\mathcal T(\mathbb K_{m,n})$ is
$$
\sum_{\substack{p,a,b\ge 0\\ p+a+b=d}}\binom{m}{a}\binom{n}{b}\,(p!)^2\,S(m-a,p)\,S(n-b,p)\,c_{p+1}(a,b).
$$
In particular, the number of regions is $c_1(m,n)$, the number of acyclic orientations of $\mathbb K_{m,n}$.
\end{theorem}

\begin{proof}
To count the number of $d$-dimensional faces of $\mathcal T(\mathbb K_{m,n})$, we will count the number of colored $\mathbb K_{m,n}$-threshold graphs with $d$ color classes. By the proof of Lemma~\ref{lem:orient}, every color class is a single vertex or induces a complete bipartite graph, so it is bipartite. By Theorem~\ref{thm:dim_through_thresh_grph}, the dimension of the face is equal to the number of color classes.

Let $a$ be the number of color classes with one vertex in $X$, $b$ the number with one vertex in $Y$, and $p$ the number with at least two vertices, so $a+b+p=d$. The set of color classes can be chosen in
$$
\binom{m}{a}\binom{n}{b}\,p!\,S(m-a,p)\,S(n-b,p)
$$
ways: we choose the singletons, we split the remaining vertices of $X$ and of $Y$ into $p$ blocks each, and we match the blocks of $X$ with the blocks of $Y$.

By Lemma~\ref{lem:orient}, for a fixed set of color classes the colored graphs are counted by the acyclic orientations of $\mathbb{K}_p\vee\mathbb K_{a,b}$. By Stanley's theorem~\cite{Stanley1973}, the number of acyclic orientations of a graph $H$ is $\lvert\chi_H(-1)\rvert$, where $\chi_H$ is the chromatic polynomial of $H$. In a proper coloring of $\mathbb{K}_p\vee\mathbb K_{a,b}$ with $q\geq p$ colors, the vertices of $\mathbb{K}_p$ receive distinct colors, and since every vertex of $\mathbb K_{a,b}$ is adjacent to every vertex of $\mathbb{K}_p$, the graph $\mathbb K_{a,b}$ is colored properly with the remaining $q-p$ colors. Hence
$$
\chi_{\mathbb{K}_p\vee\mathbb K_{a,b}}(q)=q(q-1)\cdots(q-p+1)\,\chi_{\mathbb K_{a,b}}(q-p),
$$
first for integers $q\ge p$, and then as an identity of polynomials. To compute $\chi_{\mathbb K_{a,b}}$, we color the $a$ vertices of $X$ with exactly $i$ colors, which can be done in $S(a,i)\,q(q-1)\cdots(q-i+1)$ ways (split $X$ into $i$ blocks, then give the blocks distinct colors), and then each of the $b$ vertices of $Y$ with any of the other $q-i$ colors. Thus
$$
\chi_{\mathbb K_{a,b}}(q)=\sum_{i=0}^{a}S(a,i)\,q(q-1)\cdots(q-i+1)\,(q-i)^{b}.
$$
At $q=-1$, the first factor has absolute value $p!$. Put $s=p+1$, so that $q-p=-s$. Then
$$
(-s)(-s-1)\cdots(-s-i+1)=(-1)^i s^{\overline{i}},\qquad (-s-i)^b=(-1)^b(s+i)^b.
$$
Since $(-1)^{a+b}\chi_{\mathbb K_{a,b}}(-s)>0$, we get $|\chi_{\mathbb K_{a,b}}(-s)|=c_s(a,b)$. Thus the number of acyclic orientations of $\mathbb{K}_p\vee\mathbb K_{a,b}$ is $p!\,c_{p+1}(a,b)$.

Multiplying the number of sets of color classes by the number of acyclic orientations, we obtain
$$
\binom{m}{a}\binom{n}{b}\,(p!)^2\,S(m-a,p)\,S(n-b,p)\,c_{p+1}(a,b)
$$
colored $\mathbb K_{m,n}$-threshold graphs for fixed $(p,a,b)$. Summing over all $p,a,b$ with $p+a+b=d$ gives the formula.
\end{proof}

\section*{AI disclosure}
The combinatorial model (Sections~2 and~3), the results on the threshold arrangement $\mathcal T(\mathbb K_n)$, and the overall ideas of the paper are due to the author. Claude (Anthropic) carried out the application of the model to the complete bipartite graph $\mathbb K_{m,n}$ (Section~4.2), including the correspondence with acyclic orientations and the first version of the proof of Theorem~\ref{thm:faces_Kmn}. Claude was also used to review the proofs, to suggest corrections and clarifications, and to write and run computer checks of the formulas for small cases. The free version of ChatGPT (OpenAI) was used to improve the writing. The author checked all statements and proofs and are responsible for the content of the paper.

\section*{Acknowledgements}
The main problem of this project was shared during the NCM Workshop on Topics in Geometric Combinatorics, held at CMI, Chennai, in July 2025. The author is grateful to the organizers for providing a productive and inspiring research environment, and would especially like to thank Dr.\ Priyavrat Deshpande for his valuable discussions and insightful suggestions. The author wishes to thank the Council of Scientific and Industrial Research (CSIR), New Delhi, for the award of a Junior Research Fellowship [File No.~09/0079(19932)/2024-EMR-I].
\bibliographystyle{alpha}
\bibliography{main}

@article{DeshpandeMenonSingh2021,
  title={A combinatorial statistic for labeled threshold graphs},
  author={Deshpande, Priyavrat and Menon, Krishna and Singh, Anurag},
  journal={arXiv preprint arXiv:2103.03865},
  year={2021}
}

@article{Spiro2019,
  title={Counting labeled threshold graphs with Eulerian numbers},
  author={Spiro, Sam},
  journal={arXiv preprint arXiv:1909.06518},
  year={2019}
}

@book{MahadevPeled1995,
  title={Threshold graphs and related topics},
  author={Mahadev, Nadimpalli VR and Peled, Uri N},
  volume={56},
  year={1995},
  publisher={Elsevier}
}

@book{Shi1986,
    AUTHOR = {Shi, Jian Yi},
     TITLE = {The {K}azhdan-{L}usztig cells in certain affine {W}eyl groups},
    SERIES = {Lecture Notes in Mathematics},
    VOLUME = {1179},
 PUBLISHER = {Springer-Verlag, Berlin},
      YEAR = {1986},
     PAGES = {x+307},
      ISBN = {3-540-16439-1},
   MRCLASS = {20G05 (17B10 20B25 20H15)},
  MRNUMBER = {835214},
MRREVIEWER = {H.\ H.\ Andersen},
       DOI = {10.1007/BFb0074968},
       URL = {https://doi.org/10.1007/BFb0074968},
}

@article{Athanasiadis1999,
  title={A simple bijection for the regions of the Shi arrangement of hyperplanes},
  author={Athanasiadis, Christos A and Linusson, Svante},
  journal={Discrete mathematics},
  volume={204},
  number={1-3},
  pages={27--39},
  year={1999},
  publisher={Elsevier}
}

@article{Athanasiadis1996,
  title={Characteristic polynomials of subspace arrangements and finite fields},
  author={Athanasiadis, Christos A},
  journal={Advances in mathematics},
  volume={122},
  number={2},
  pages={193--233},
  year={1996},
  publisher={Elsevier}
}

@article{Stanley2007,
  title={An introduction to hyperplane arrangements},
  author={Stanley, Richard P},
  journal={Geometric combinatorics},
  volume={13},
  pages={389--496},
  year={2007}
}

@article{Stanley1973,
  title={Acyclic orientations of graphs},
  author={Stanley, Richard P},
  journal={Discrete Mathematics},
  volume={5},
  number={2},
  pages={171--178},
  year={1973},
  publisher={Elsevier}
}

@book{Stanley2011,
    AUTHOR = {Stanley, Richard P.},
     TITLE = {Enumerative combinatorics. {V}olume 1},
    SERIES = {Cambridge Studies in Advanced Mathematics},
    VOLUME = {49},
   EDITION = {Second},
 PUBLISHER = {Cambridge University Press, Cambridge},
      YEAR = {2012},
     PAGES = {xiv+626},
      ISBN = {978-1-107-60262-5},
   MRCLASS = {05-02 (05A15 06-02)},
  MRNUMBER = {2868112},
}
\end{document}